# On Certain Bounds for Multiset Dimensions of Zero-Divisor Graphs Associated with Rings

**Nasir Ali[a], Hafiz Muhammad Afzal Siddiqui[a], Muhammad Imran Qureshi[b]**
[a]Department of Mathematics, COMSATS University Islamabad, Lahore Campus, Pakistan.
[b]Department of Mathematics, COMSATS University Islamabad, Vehari Campus, Pakistan.
**Email address:**
nasirzawar@gmail.com (Corresponding Author: Nasir Ali), hmasiddiqui@gmail.com (Hafiz Muhammad Afzal Siddiqui), imranqureshi18@gmail.com,( Muhammad Imran Qureshi)

**Abstract:**

This article investigates multiset dimensions (Mdim) in zero divisor graphs (ZD-graphs) associated with rings. Through rigorous analysis, we establish general bounds for the multiset dimension in ZD-graphs, exploring various commutative rings including the ring of Gaussian integers, the ring $Z_n$ of integers *modulo n*, and quotient polynomial rings. Additionally, we examine the behavior of Mdim under algebraic operations and discuss bounds in terms of diameter and maximum degree. By emphasizing the correlation between Mdim and the concept of symmetry, this study not only enhances our understanding of algebraic structures and their graphical representations but also highlights the intrinsic symmetrical nature inherent in these mathematical constructs.

1. Introduction:

The interplay between graph theory and algebra has led to significant insights, particularly through the study of ZD-graphs associated with commutative rings. Beck [1] pioneered this connection by introducing the concept of ZD-graphs denoted by $Z^o(R)$, focusing on the correspondence between ring elements and graph nodes, where a zero vertex is connected to all other remaining vertices. Subsequently, Anderson and Livingston [2] extended this notion to consider ZD-graphs where each vertex stands for a nonzero zero divisor. In this formulation, an undirected graph is formed by considering ring elements $x$ and $y$ as nodes connected by an edge if their product equals zero. Notably, Anderson and Livingston's investigation concentrated on finite rings, establishing connections between ring properties and graph theoretic properties such as completeness or star structure. This approach, denoted as $Z(R)$, slightly deviates from Beck's original definition of ZD-graphs, wherein zero is not considered a vertex. Anderson and Livingston's findings shed light on the relationship between ring properties and the structural characteristics of $Z(R)$, offering valuable insights into both algebraic and graph theoretic domains.

Expanding upon this foundation, Redmond [3] extended the study of ZD-graphs to include noncommutative rings, introducing various methods to characterize ZD-graphs for both undirected and directed graphs. Redmond [4] further enriched this exploration by introducing ideal-based ZD-graphs for commutative rings. This novel approach involves substituting elements with zero products with elements belonging to a specific ideal I of the ring R, yielding a graph denoted as $\Gamma_I(R)$. Throughout the paper $L(R)$ denotes the set of zero divisors of ring $R$ and $Z(R)$ will be type of ZD-graph under consideration.

Throughout the literature, numerous variations of ZD-graphs have been proposed, including total graphs, unit graphs, Jacobson graphs, and ZD-graphs based on equivalence classes. Each variant offers a distinct perspective on the algebraic and structural properties inherent in commutative rings, enriching the study of ZD-graphs and their applications in both algebra and graph theory. These works can be found in sources such as [5-9]. Readers may refer to [10, 11] for a deeper understanding of graph theory, and for some fundamental definitions of ring theory, [12, 13] can be consulted.

Simanjuntak et al. [15] introduced a new variant of metric dimension (md) known as multiset dimension

(*Mdim*), where a multiset of distances between $v$ and all vertices in resolving set (RS) $W$ were calculated including their multiplicities. The *Mdim* is defined as the minimum cardinality of the RS. Assume $v$ is a node of $G$ and $B \subseteq V(G)$, then representation for multiset of $v$ w.r.t $B$ defined as the distances between $v$ and nodes in $B$. This representation is denoted by $r_m(B)$. For every pair of distinct nodes $u$ and $v$, $B$ is called $m$ −resolvingset (M-RS) of $G$, if $r_m(B) \neq r_m(B)$. The cardinality of M-RS is multiset basis of $G$ and minimum cardinality of multiset basis is called multiset dimension of $G$, denoted by $Mdim(G)$, if $G$ does not contain a M-RS, we write $Mdim(G) = \infty$. The key point of this article is that the apparent expansions are an oversimplification of the task of identifying graph vertices using the multiset representation.

The article presents significant contributions to the field. Primarily, it expands the notion of Mdim to ZD-graphs, with a special emphasis on the ZD-graph representing the ring of integers $Z_n$ modulo $n$. This expansion facilitates varied descriptions of rings by examining the vertices within their respective ZD-graphs, integrating metric and multiset representations effectively. Through rigorous proofs, it establishes that rings can be effectively characterized by their multiset dimensions. Additionally, it offers insights into various rings, demonstrating that their ZD-graphs can have bounded multiset dimensions determined by the graph's diameter. Moreover, the article presents a straightforward method for computing the Mdim of ZD-graphs for rings of integers modulo n, enhancing accessibility to this important metric. The novelty of this research lies in its exploration of multiset dimensions in the context of graph theory, a relatively underexplored area. By unveiling the multiset dimensions of ZD-graphs, the research significantly enhances our grasp of the inherent structural and algebraic characteristics found within these graphs. This progress bodes well for various practical applications such as network design, social networking, and communication systems. A deeper comprehension of graph structures is imperative for maximizing performance and efficiency in these domains.

## 2. Preliminaries

In the preliminary section of this research article, we establish foundational concepts and terminology pertinent to graph theory and its applications within algebraic structures. A graph, formally defined as an ordered pair $G = (V, E)$, comprises a set of vertices or nodes ($V$) and a set of edges ($E$). The order of a graph refers to the cardinality of its vertex set, while its size pertains to the cardinality of its edge set. The distance between two nodes $u'$ and $v'$ is defined as the length of the shortest path connecting them, while the distance from a node $w$ to an edge $e' = u'v'$ is determined as the minimum of the distances from $w$ to $u'$ and $v'$. Graphs exhibit diverse structural properties, including regularity and completeness. A graph is considered regular if every vertex has the same degree, specified as $deg(r) = c$ for a fixed $c$ in the positive integers. Complete graphs, denoted as $k_m$, establish connections between all pairs of vertices, where $m$ denotes the number of vertices. Complete bipartite graphs, typically represented as $k_{m,n}$, are partitioned into two distinct sets of vertices, $X$ and $Y$, with each vertex in $X$ connected to every vertex in $Y$. A significant graph theoretic concept is that of a cut vertex, which arises when the removal of a vertex from a connected graph results in two or more disconnected components.

Recent research has extended graph theory into algebraic realms, particularly concerning ZD-graphs. Redmond's investigations in [17] and [18] explored ZD-graphs of noncommutative and commutative rings, respectively. Subsequent studies, including those by Siddiqui et al. in [20] and Pirzada and Aijaz in [21], have further delved into metric parameters for ZD-graphs, examining dimensions and bounds associated with these algebraic structures. Such endeavors not only enrich the theoretical understanding of graph theory but also offer practical insights into the structural properties of algebraic systems. Simanjuntak et al. [15] found some sharp bounds for *Mdim* of arbitrary graphs in term of its md, order, or diameter. Siamanjuntak also provided some necessary conditions for a graph to have finite *Mdim*, with an example of an infinite family of graphs where those necessary conditions are also sufficient. It was also shown that the *Mdim* of any graph other than a path is at least 3 and two families of graphs having the *Mdim* 3 were

proved. Here, we consider some results from [15] as follows:

**Theorem 2.1 [15]:** For any integer $m \geq 3, n \geq 6$, $Mdim(P_m) = 1$, $Mdim(k_m) = \infty$ and $Mdim(C_n) = 3$. Moreover, $Mdim(P_m) = 1$ iff $G \cong P_m$.

Moreover, for a complete bipartite graph $k_{m,n}$ we get different $Mdim$ for different choices of values of $m$ & $n$.

**Theorem 2.2 [15]:** For any complete bipartite graph $k_{m,n}$, the $Mdim(k_{m,n})$ is given below.

| $Mdim(k_{m,n}) = 1$ | for $m = 1$ & $n = 1,2$ |
|---|---|
| | for $m = 2$ & $n = 1$ |
| $Mdim(k_{m,n}) = \infty$ | for $m = 1$ & $n \geq 3$ |
| | for $m = 2$ & $n \geq 2$ |

Moreover, $Mdim$ for a single vertex graph $G$ is supposed to be zero and for an empty graph it is undefined. Our discussion commences with the subsequent observation.

## 3. Results and Discussion

**Theorem 3.1:** Let $R$ be a finite commutative ring with unity. If $R$ is an integral domain (ID), then $\boldsymbol{Mdim(Z(R))}$ is undefined.

*Proof:*

It is well-known that if R is an integral domain (ID), then the set of zero divisors $Z(R)$ is not defined, implying that the multiset dimension (Mdim) of the ZD-graph $Z(R)$ is undefined, and conversely. Now, the following result provides the Mdim of the ZD-graph of ring $R$ when the set $Z(R)$ is isomorphic to the path graph $P_m$. ∎

**Proposition 3.1:** Consider $R$ as a finite commutative ring having unity. Then $Mdim(Z(R)) = 1$ iff $R$ is isomorphic to one of the following rings: $\mathbb{Z}_6, \mathbb{Z}_8, \mathbb{Z}_9, \mathbb{Z}_2 \times \mathbb{Z}_2, \mathbb{Z}_3(r)/(r^2), \mathbb{Z}_2(r)/(r^3)$, or $\mathbb{Z}_4(r)/(2r, r^2 - 2)$.

*Proof:* Suppose that $Mdim(Z(R)) = 1$. The ZD-graphs of rings $\mathbb{Z}_6, \mathbb{Z}_8, \mathbb{Z}_9, \mathbb{Z}_2 \times \mathbb{Z}_2, \mathbb{Z}_3(r)/(r^2), \mathbb{Z}_2(r)/(r^3)$, or $\mathbb{Z}_4(r)/(2r, r^2 - 2)$ are path graphs, and we know that path graphs are only graphs whose $Mdim = 1$ by Theorem 2.1. Since $|L(R)|$ is at most three whenever $Z(R) \cong P_n$. Using ([19], Lemma 2.6), $Z(R)$ is either $P_2$ or $P_3$.

Case I: Let's consider the case where $Z(R) \cong P_2$, moreover, let $|L(R)| = a, b$ satisfying $a.b = 0$. Rings that fulfill this property encompass $\mathbb{Z}_9, \mathbb{Z}_2 \times \mathbb{Z}_2, \mathbb{Z}_3(r)/(r^2)$.

Case II: In the case where $Z(R)$ is isomorphic to $P_3$, considering $|L(R)| = a, b, c$, satisfying $a.b = 0$ and $b.c = 0$. Rings satisfying these conditions include $\mathbb{Z}_6, \mathbb{Z}_8, \mathbb{Z}_2(r)/(r^3)$, and $\mathbb{Z}_4(r)/(2r, r^2 - 2)$ [14].

Conversely, the ZD-graphs of the mentioned rings are either $P_2$ or $P_3$, [14]. Thus, by Theorem 2.1, $Mdim(Z(R)) = 1$. ∎

**Proposition 3.2:** Let $R$ be a finite commutative ring with unity and $R$ is isomorphic to one of the following rings, $\mathbb{Z}_3 \times \mathbb{Z}_3$, $K_4(r)/(r^2), \mathbb{Z}_4(r)/(r^2+r+1), \mathbb{Z}_4(r)/(2,r)^2, \mathbb{Z}_2[r,s]/(r,s)^2$. Then $Z(R) \cong C_m$, and $Mdim(Z(R))$ is infinite.

*Proof:* Assuming $R$ is a commutative ring with unity and $Z(R)$ is a cyclic graph, according to ([3], Theorem 2.4), it follows that the maximum length of the cyclic graph is 4. Since $Z(R) \cong C_m$ and $m$ does not exceed 4, so by Theorem 2.1, $Mdim(Z(R)) = \infty$.

The corresponding ZD-graphs for the rings mentioned above can be observed in Figure 1. ∎

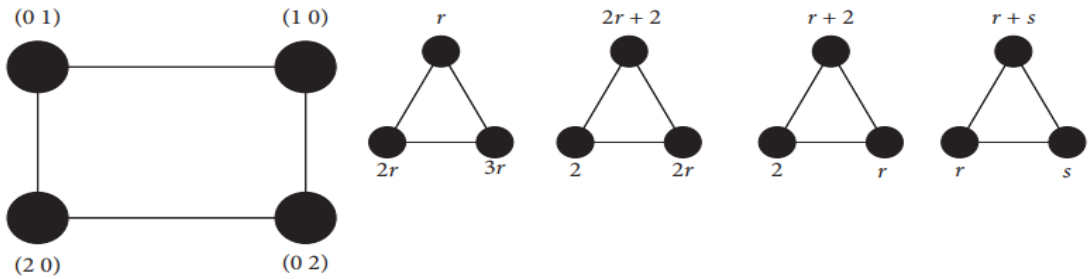

*Figure 1: ZD-graphs of* $\mathbb{Z}_3 \times \mathbb{Z}_3$, $K_4(r)/(r^2), \mathbb{Z}_4(r)/(r^2+r+1), \mathbb{Z}_4(r)/(2,r)^2, \mathbb{Z}_2[r,s]/(r,s)^2$

**Theorem 3.2:** Let $R$ be a finite commutative ring with unity. If $R$ is isomorphic to one of the following rings $\mathbb{Z}_2 \times \mathbb{Z}_4$, $\mathbb{Z}_2 \times \mathbb{Z}_4[X]/(x^2)$, $\mathbb{Z}_2 \times \mathbb{Z}_2 \times \mathbb{Z}_2$, then $Mdim(Z(R)) = 3$.

*Proof:* Consider the ring $R$ is isomorphic to one of the following rings $\mathbb{Z}_2 \times \mathbb{Z}_4$, $(\mathbb{Z}_2 \times \mathbb{Z}_4)[X]/(x^2)$, $\mathbb{Z}_2 \times \mathbb{Z}_2 \times \mathbb{Z}_2$. The ZD-graphs of the above rings are shown in Fig. 3 and Fig. 4. then there exists a minimal resolving set for $Z(R)$ say, $\{v_1, v_2, v_3\}$ by ([3], Theorem 2.3). Now by ([23], Theorem 2.3) $Mdim(Z(R)) = 3$. ∎

**Theorem 3.3:** Let R be a finite commutative ring with unity, such that each element $r$ in the zero-divisor set $L(R)$ is nilpotent. Then:

a) If $|L(R)| \geq 3$ and $L(R)^2 = \{0\}$, then $Mdim(Z(R)) = \infty$.

b) If $|L(R)| \geq 3$ and $L(R)^2 \neq 0$, then $Mdim(Z(R))$ is finite.

*Proof:* a) Assuming that $|L(R)| \geq 3$ and $L(R)^2 = \{0\}$, it follows that the product of any pair of elements $r$ and $s$ in L(R) is zero, i.e., $r.s = 0$ for all $r, s \in L(R)$ and By [5], Theorem 2.8, this condition implies that $Z(R)$ forms a complete graph. Therefore, according to Theorem 2.1, the Mdim of $Z(R)$ is infinite.

b) Considering $|L(R)|$ is greater than or equal to 3 and $L(R)^2 \neq 0$ there exists an element $r$ in $L(R)$ such that $r^2 = 0$, indicating the occurrence of some element $s$ in $L(R)$ thereby the distance between $r$ and $s$ is at least 2. Consequently, $L(R)/(r,s)$ serves as a multiset generator for any vertex $s$ adjacent to $r$. Thus, the Mdim of $Z(R)$ is finite in this case.

**Corollary 3.1:** Let R be a finite commutative ring with unity such that $|L(R)| \geq 3$. If Z(R) has a cut vertex but no degree 1 vertex, then $Mdim(Z(R)) = \infty$.

*Proof:* For the given ring $R$, suppose the ZD-graph $Z(R)$ has a cut vertex but no vertex with degree 1. According to ([25], Theorem 3), this implies that $R$ exhibits an isomorphism with one of the following rings: $\frac{\mathbb{Z}_2[r,s]}{(r^2,s^2-rs)}, \frac{\mathbb{Z}_4[r]}{(r^2+2r)}, \frac{\mathbb{Z}_4[r,s]}{(r^2,s^2-rs,\ rs-2,\ 2r,\ 2s)}, \frac{\mathbb{Z}_8[r,s]}{(2r,\ r^2+4)}, \frac{\mathbb{Z}_2[r,s]}{(r^2,s^2)}, \frac{\mathbb{Z}_4[r]}{(r^2)}, \frac{\mathbb{Z}_4[r,s]}{(r^2,s^2,\ rs-2,2r,2s)}$.

One should note that Theorem 3 in [25] establishes that if the ZD-graph has a cut vertex but no vertex with degree one, the authors provide a list of rings, have associated ZD-graphs following this property.

Figure 2(a) which represents the $Z(R)$ for the first four rings, in Figure 2(b), representing $Z(R)$ for the remaining three rings. Therefore, it is concluded that $Mdim(Z(R)) = \infty$. ∎

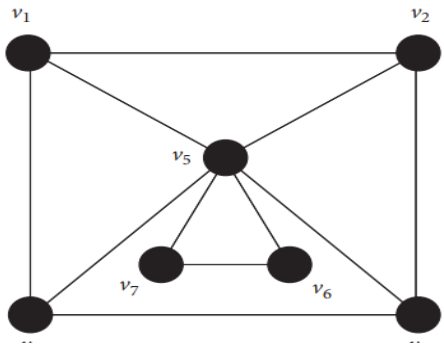
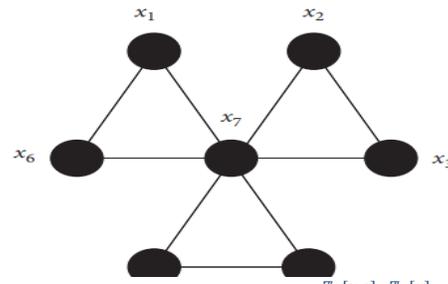

Figure 2 (b): ZD-graph of Rings $\frac{\mathbb{Z}_2[r,s]}{(r^2,s^2)}, \frac{\mathbb{Z}_4[r]}{(r^2)}, \frac{\mathbb{Z}_4[r,s]}{(r^2,s^2,\ rs-2,2r,2s)}$

Figure 2 (a): ZD-graph of Rings $\frac{\mathbb{Z}_2[r,s]}{(r^2,s^2-rs)}, \frac{\mathbb{Z}_4[r]}{(r^2+2r)}, \frac{\mathbb{Z}_4[r,s]}{(r^2,s^2-rs,\ rs-2,\ 2r,\ 2s)}, \frac{\mathbb{Z}_8[r,s]}{(2r,\ r^2+4)}$

**Corollary 3.2:** Let $R$ be a finite commutative ring with unity and $R$ be a local ring having no cycles in associated $Z(R)$ then $Mdim(Z(R)) = 1$.

**Proof.** Let $R$ is a local ring having no cycle in associated Z(R), then by ([25], Theorem 2.1), Z(R) is isomorphic to $P_2$ or $P_3$. Hence $Mdim(Z(R)) = 1$. ∎

**Theorem 3.5:** Let $\mathbb{Z}_n$ be a ring of integers modulo $n$. Assuming $p$ and $q$ as distinct primes, we have $Mdim(Z(\mathbb{Z}_n))$ as follows:

| $n$ | $Mdim(Z(\mathbb{Z}_n))$ |
| --- | --- |
| $p$ | undefined |
| $p^2$ | $\infty$ |
| $pq$ | $\infty$ |
| $2^2$ | 0 |
| $2^3$ | 1 |
| $3^2$ | 1 |

*Proof:* To prove the theorem, it is imperative to address multiple cases individually.

(i) When $n = 2^2$, then $Z(\mathbb{Z}_n)$ consists of single vertex. So $Mdim(Z(\mathbb{Z}_n)) = 0$. Now consider $n = 3^2, 2^3$, then $Z(\mathbb{Z}_n)$ consists of two and three vertices respectively, and $Z(R) \cong P_n$ (where $n = 2$ or $3$). By Theorem 2.1, we know that $Mdim(P_m) = 1$ which follows that $Mdim(Z(\mathbb{Z}_n)) = 1$. Suppose $n = pq$, where $p$ and $q$ are distinct primes. We can partition the vertices into two sets: $U = \{\lambda p \in Z(\mathbb{Z}_n); (\lambda, q) = 1\}$ and $V = \{\lambda q \in Z(\mathbb{Z}_n); (\lambda, p) = 1\}$. This partition clearly demonstrates that $Z(\mathbb{Z}_n)$ is bipartite. Furthermore, it is evident that $x.y = 0$ for every $x \in U$ and $y \in V$. Therefore $Z(\mathbb{Z}_n)$ forms a complete bipartite graph. Consequently, when $n = pq$, $Z(\mathbb{Z}_n) \cong K_{q-1,p-1}$. Thus, by Theorem 2.2, $Mdim(Z(\mathbb{Z}_n)) = \infty$.

(ii) Consider $p > 3$, if we take $n$ different from case (i), i.e., if $n = p$, then $Z(\mathbb{Z}_n) = \varphi$ so $Mdim(Z(\mathbb{Z}_n)) = \infty$. If $n = p^2, p^n$ or $2^2p$ or any other $n$, then $Z(\mathbb{Z}_n) \cong K_{p-1}, K_{n,m}$ (*not complete*), *and* $K_{q-1,p-1}$. By Theorem 2.2, $Mdim(Z(\mathbb{Z}_n)) = \infty$. ∎

Table 1: Multiset dimension of $Z(\mathbb{Z}_n)$

| n | \|V\| | \|E\| | Diameter | Girth | $Z(\mathbb{Z}_n)$ | $Mdim(Z(\mathbb{Z}_n))$ |
|---|---|---|---|---|---|---|
| p | 0 | 0 | 0 | Undefined | $Z(\mathbb{Z}_n) = \emptyset$ | ∞ |
| $2^2$ | 1 | 0 | 0 | Undefined |  | 0 |
| $3^2$ | 2 | 1 | 1 | Undefined | 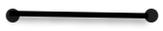 | 1 |
| $P^2$, $p \geq 5$ | p-1 | $\binom{p-1}{2}$ | 1 | 3 | Complete graph $K_{p-1}$ | ∞ |
| $2^3$ | 3 | 2 | 2 | Undefined | 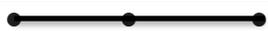 | 1 |
| $P^n$, $n \geq 3$ | $P^{n-1}-1$ | $\sum (p-1)^{\lfloor i/2 \rfloor}$ | 2 | 3 | Quasi-local: $P^{n-1}$ is attached to everything | ∞ |
| $2^2p$, $p \geq 3$ | 2p+1 | 4p-4 | 3 | 4 | Bipartite graph (but not complete) | ∞ |
| pq | q-1+p-1 | (q-1)(p-1) | 2 | 4 | Complete bipartite graph $K_{q-1,p-1}$ | 1 For $(p=3, q=2)$ ∞ for all other values |
| All others | | | 2 | 3 | | |

## 4. Bound between multiset dimensions and diameter of Zero-divisor graphs:

Within this section, we explore the relationship between Mdim and the diameter of ZD-graphs. Through our analysis, we uncover insightful bounds that enhance our understanding of these graph properties and their structural nuances.

**Lemma 4.1:** Suppose $G$ is a ZD-graph with a diameter of at most 2, and $Z(R) \not\cong P_n$, then the multiset dimension of $Z(R)$ is infinite.

*Proof:* Consider the ZD-graph $G$ with a maximum diameter of 2, and assume $Z(R)$ is not a path graph. Consequently, $Z(R)$ can take the form of a cycle with a maximum of 5 vertices, a complete graph, the Peterson graph, or a star graph. Then by using Theorem 2.1 and 2.2 and Theorem 3.1 [15] the result follows.
∎

**Corollary 4.1:** Let $Z(R)$ be a ZD-graph of a ring $R$ such that $|L(R)| \geq 3$ and diameter $d$. Then,

$$Mdim(Z(R)) > f(n, d).$$

(For positive integers $n, d$, $f(n, d)$ is defined to be least positive integer k for which $\frac{(k+d-1)!}{k!(d-1)!} + k \geq n$ )

*Proof:* Consider a commutative ring $R$ with identity and let $Z(R)$ denote its zero-divisor graph. For positive integers $n$ and $d$, the function $f(n, d)$ is defined to be the least positive integer $k$ for which,

$$\frac{(k + d - 1)!}{k!\,(d - 1)!} + k \geq n$$

We aim to prove that the multiset dimension $Mdim(Z(R))$ is at least some function $f(n, d)$.

For $n = 3$, the diameter $d$ of the zero divisor graph can be at most 2, as for any three vertices, there exists a pair with a distance of at most 2 between them. Thus, $Mdim(Z(R)) = 2$.

For $n > 3$, let's examine the definition of $f(n, d)$. We want to find the least positive integer $k$ satisfying the inequality mentioned above. This inequality essentially gives us a lower bound on the cardinality of $S$ (resolving set), which is $k$. Thus, $Mdim(Z(R)) \geq k$.

Therefore, $Mdim(Z(R))$ is at least $f(n, d)$. Alternatively, the proof is followed by using Theorem 2.3 and Theorem 2.4 [15]. ∎

**Theorem 4.2:** Consider a ZD-graph $Z(R)$ with $n$ vertices, and having diameter $k$, if $Mdim(Z(R)) = 3$, then

$$n \leq \frac{k^2(k + 3) + 2(k + 6)}{6}$$

*Proof:* We can deduct from the observation that no vertex in a ZD-graph of multiset dimension 3 can have the representation as $(0, 1, 1)$. Like this, no vertex may have the representation $(1, 1, 1)$.
The interpretation of vertex $v \in V$ with respect to resolving set $W$ outside the $W$ will take the form of

$\{1^{a_1}, 2^{a_2}, \ldots, k^{a_k}\}$, such that $a_1 + a_2 + a_k = 3$. This equation has number of solutions as follows:
$\binom{k+2}{k-1} = \frac{k^2(k+3)+2k}{6}$.

Since $(1, 1, 1)$ should not be included, then possible number of representations for vertices outside the resolving set is,
$$\frac{k^2(k+3)+2k}{6} - 1$$
As $(n-3)$ vertices are outside the resolving set. Then,
$$\frac{k^2(k+3)+2k}{6} - 1 \geq n - 3$$
Hence,
$$n \leq \frac{k^2(k+3)+2k}{6} - 1 + 3 = \frac{k^2(k+3)+2(k+6)}{6}. \blacksquare$$

Below are some of the examples of the ZD-graphs of rings having multiset dimensions 3.

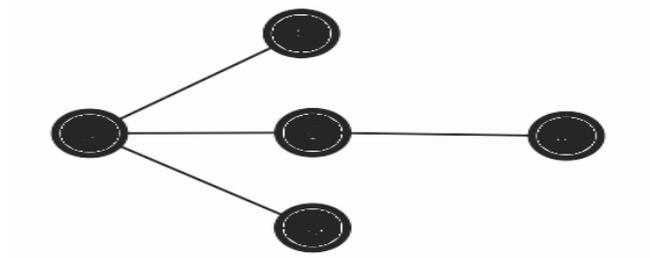

Figure 3: Graph of $\mathbb{Z}_2 \times \mathbb{Z}_4$ and $\mathbb{Z}_2 \times \mathbb{Z}_4[X]/(x^2)$

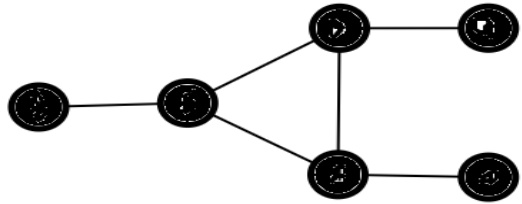

Figure 4: Graph of $\mathbb{Z}_2 \times \mathbb{Z}_2 \times \mathbb{Z}_2$

**Open Problem:** Find a family of ZD-graphs of rings which satisfy the equality in Theorem 4.2?

5. **Methodology**

To investigate the multiset dimensions of ZD-graphs associated with rings, we employed a comprehensive methodology involving data visualization, categorization, and dimension analysis. Initially, we utilized both MATLAB and Python algorithms to visualize the ZD-graphs corresponding to various commutative rings, including the ring $Z_n$ of integers modulo $n$, the ring of Gaussian integers modulo $m$, and quotient polynomial rings. This visualization process allowed us to gain insights into the structural properties of the graphs and identify key characteristics. Subsequently, we undertook a rigorous categorization process to classify the ZD-graphs based on their structural features. This categorization facilitated the systematic analysis of multiset dimensions across different types of ZD-graphs, enabling us to identify patterns and variations among them.

Following the categorization phase, we conducted in-depth dimension analysis to determine the multiset dimensions of the categorized ZD-graphs. This involved implementing algorithms to calculate the multiset dimensions and establish relationships between the dimensions and other graph parameters, such as diameter and maximum degree. Throughout the methodology, meticulous attention was paid to ensure the accuracy and reliability of the analysis results. The use of both MATLAB and Python algorithms provided flexibility and robustness in data processing and visualization. Additionally, the categorization and dimension analysis procedures were carried out systematically to ensure comprehensive coverage and thorough exploration of the research objectives.

Overall, the methodology employed in this study facilitated a detailed investigation into the multiset dimensions of ZD-graphs associated with rings, enabling us to uncover valuable insights into the algebraic and graph-theoretic properties of these structures.

## 6. Discussion

The exploration of Mdim in ZD-graphs presents a compelling avenue for understanding the structural properties of rings and their graphical representations. In this section, we delve into the implications of our findings and discuss the potential applications and future directions of this study. Our analysis reveals that the multiset dimension serves as a distinctive feature that can aid in characterizing rings based on their associated zero divisor graphs. By establishing general bounds and exploring the behavior of Mdim across different types of rings, we can discern characteristic patterns and attributes. This classification approach offers a systematic method for identifying and distinguishing rings based on their algebraic structures. The ability to characterize rings through multiset dimensions has practical implications in various domains. In algebraic cryptography, for instance, understanding the Mdim of rings can inform the selection of cryptographic protocols, ensuring the security and efficiency of cryptographic systems. Additionally, in computational algebra, knowledge of the structural properties of rings via ZD-graphs can streamline algorithms for solving algebraic equations and optimizing computational processes. By examining the relationship between multiset dimensions and zero divisor graphs, we deepen our understanding of algebraic structures and their graphical representations. This study sheds light on the intricate interplay between algebraic concepts and graph theory, highlighting the rich connections between seemingly disparate mathematical disciplines. Moreover, the exploration of Mdim provides insights into the inherent properties and behaviors of rings, paving the way for further research and inquiry. Looking ahead, there are several avenues for future research in this area. Expanding our analysis to encompass a broader range of rings and exploring additional properties of zero divisor graphs could yield further insights into the relationship between algebraic structures and graphical representations. Moreover, investigating the applicability of Mdim in other contexts, such as network theory or combinatorics, could uncover new connections and applications beyond the realm of algebra. In summary, the study of multiset dimensions in zero divisor graphs offers a fruitful avenue for exploring the structural properties of rings and their graphical representations. By leveraging Mdim as a characterization tool, we can identify rings based on their distinctive features and discern characteristic patterns across different types of rings. This discussion not only enhances our understanding of algebraic structures but also opens up new avenues for practical applications and future research endeavors.

## 7. Conclusion

In this study, we have investigated the characterization of rings through the analysis of multiset dimensions associated with their ZD-graphs. Specifically, we have explored the multiset dimensions of ZD-graphs corresponding to various commutative rings, such as the ring $Z_n$ of integers modulo $n$, the ring of Gaussian integers modulo $m$, and polynomial rings. Additionally, we have generalized the multiset dimension for the ring $Z_n$ of integers modulo $n$. Our investigation has culminated in the presentation of bounds relating the multiset dimension to the diameter of ZD-graphs. Moving forward, there are several avenues for future research in this domain. Firstly, extending the analysis to include noncommutative rings could offer valuable insights into the structural properties of ZD-graphs in a broader algebraic context. Moreover, exploring the multiset dimensions of ZD-graphs associated with other algebraic structures beyond rings could further enrich our understanding of their graph-theoretic properties. Additionally, investigating the relationships between multiset dimensions and other graph parameters, such as connectivity or chromatic numbers, could provide deeper insights into the interplay between algebra and graph theory in diverse settings. These future directions hold promise for advancing both theoretical understanding and practical applications in algebraic graph theory.

**Declaration:**

- **Availability of data and materials:** The data is provided on the request to the authors.
- **Conflicts of interest:** The authors declare that they have no known competing financial interests or personal relationships that could have appeared to influence the work reported in this paper.
- **Fundings:** No funding received.
- **Author's contribution:** All the authors equally contributed towards this work.